\documentclass[12pt]{article}

\usepackage{amssymb}
\usepackage{amsmath,amsthm}
\usepackage[dvips]{graphicx}
\usepackage{wrapfig}
\usepackage{bm}

\newtheorem{theorem}{Theorem}[section]
\newtheorem{corollary}[theorem]{Corollary}
\newtheorem{lemma}[theorem]{Lemma}
\newtheorem{proposition}[theorem]{Proposition}

\theoremstyle{definition}

\theoremstyle{remark}
\newtheorem{remark}[theorem]{Remark.}

\title{Axis systems of link projections}

\author{
Kirara Horiguchi \thanks{Department of Civil Engineering, National Institute of Technology, Gunma College, 580 Toriba-cho, Maebashi-shi, Gunma 371-8530, Japan. }
\and 
Ayaka Shimizu \thanks{Department of Mathematics, National Institute of Technology, Gunma College, 580 Toriba-cho, Maebashi-shi, Gunma 371-8530, Japan. Email: shimizu@nat.gunma-ct.ac.jp }
\and 
Ryohei Watanabe \thanks{Department of Information Science, University of Fukui, 3-9-1, Bunkyo, Fukui-shi, Fukui, 910-8507, Japan. }
\and 
Yoshiro Yaguchi \thanks{Department of Mathematics, National Institute of Technology, Gunma College, 580 Toriba-cho, Maebashi-shi, Gunma 371-8530, Japan. Email: yaguchi-y@nat.gunma-ct.ac.jp }
}
\date{\today}

\begin{document}

\maketitle

\begin{abstract}
An axis of a link projection is a closed curve which lies symmetrically on each region of the link projection. 
In this paper we define axis systems of link projections and characterize axis systems of the standard projections of twist knots. 
\end{abstract}

\section{Introduction}

\begin{figure}[ht]
\begin{center}
\includegraphics[width=80mm]{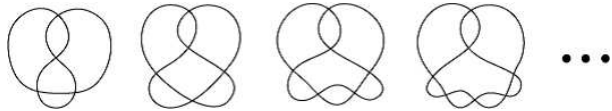}
\caption{The standard projections of twist knots. }
\label{twist}
\end{center}
\end{figure}

A {\it link} is an embedding of some circles in $S^3$, and a {\it link projection} is an image of a projection of a link on $S^2$ with double points, where arcs intersect transversely. 
We call such a double point a {\it crossing}, and each segment of a link projection between crossings which has no crossings in the interior an {\it edge}. 
We call each part of $S^2$ bounded by a link projection a {\it region}. 
We also call a link with just one component a {\it knot} and an image of its projection a {\it knot projection} (or a {\it spherical curve}). 
In this paper, we consider link projections up to ambient isotopy in $S^2$ and reflections, and we do not consider split link projections.
\begin{figure}[ht]
\begin{center}
\includegraphics[width=120mm]{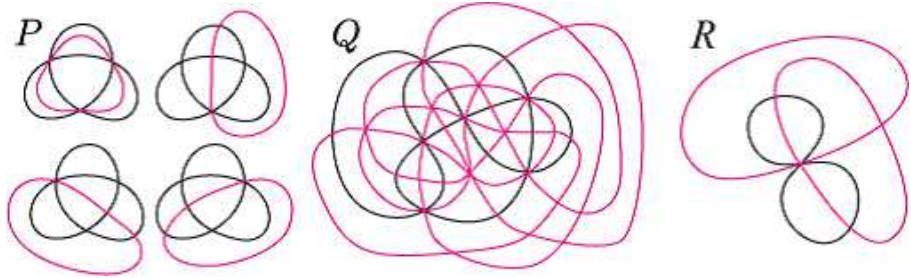}
\caption{Axes of knot projections $P$, $Q$ and $R$. }
\label{ex-axis}
\end{center}
\end{figure}
An {\it axis} of a link projection $P$, defined in Section 2, is a closed curve on $S^2$ which lies symmetrically on each region of $P$ (see Fig. \ref{ex-axis}). 
As we discuss in Appendix A, the number of axes is an invariant of link projections. 
By giving a letter to each region of $P$ and reading up the letters according to an axis starting from any point with any orientation, and giving signs, we obtain a sequence. 
We call the set $s(P)$ of the sequences (up to cyclic permutations and a reversing) for all the axes of $P$ the {\it axis system of $P$}. 
Precise definition of axis system is given in Section 3. 
The {\it standard projection of a twist knot} is the knot projection illustrated in Fig. \ref{twist}.
Note that twist knots are the knots named $4_1$, $5_2$, $6_1$, $7_2$, $8_1, \dots$.
In this paper, we characterize axis systems of the standard projections of twist knots. 

\phantom{x}
\begin{theorem}
A link projection $P$ is the standard projection of a twist knot with $n$ crossings if and only if 
\begin{center}
\begin{align*}
\bullet  s(P)= \{ & -C_0 R_0 R_1 R_2 \dots R_{m-3} R_{m-2} C_1 L_{m-2} L_{m-3} \dots  L_2 L_1 L_0, \\
& C_0 I C_1 O, -L_0 R_0 I O R_0 L_0 O I, -R_1 IO, -R_2 IO, \dots , -R_{m-2} IO, \\
& -L_1 IO, -L_2 IO, \dots , -L_{m-2} IO \}  \ ( n=2m, m=2,3,4, \dots ), \\
\bullet  s(P)= \{ & -C_0 R_0 R_1 R_2 \dots R_{m-2} R_{m-1} L_{m-1} L_{m-2} \dots L_2 L_1 L_0, \\
& C_0 I O, -L_0 R_0 I L_{m-1} O R_0 L_0 O R_{m-1} I, IO, IO, \dots ,IO, \\
& -L_1 I R_{m-2} O, -L_2 I R_{m-3} O, -L_3 I R_{m-4} O, \dots , -L_{m-2} I R_1 O \}  \\ 
& ( n=2m+1, m=2,3,4, \dots ),
\end{align*}
\end{center}
where for the case of $n=2m+1$, there are $m-1$ ``$IO$''s in $s(P)$. 
\label{mainthm}
\end{theorem}

\phantom{x}

The rest of this paper is organized as follows: 
In Section 2, we define an axis of a link projection. 
In Section 3, we define the axis system of a link projection and prove Theorem \ref{mainthm}. 
In Appendixes, we also define axis polynomials, and using them we discuss symmetries of knot projections.

\section{Axes of link projections}

\begin{figure}[ht]
\begin{center}
\includegraphics[width=40mm]{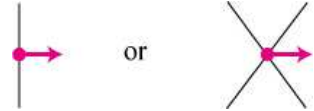}
\caption{A start point at an edge or a crossing. }
\label{startpt}
\end{center}
\end{figure}

\begin{figure}[ht]
\begin{center}
\includegraphics[width=90mm]{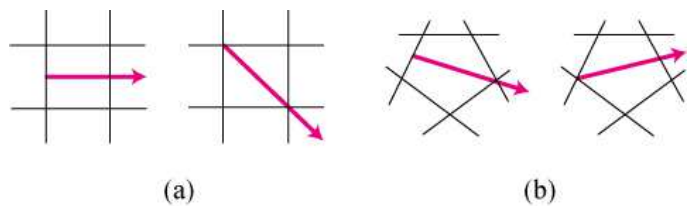}
\caption{Travel symmetrically. }
\label{evenoddregions}
\end{center}
\end{figure}

In this section we define an axis of a link projection, and discuss its properties. 
For a region of a link projection, we call it an {\it $n$-gon} if the region has $n$ edges on the boundary. 
There are many studies about $n$-gons in knot theory (see, for example, \cite{AST, aida, shinjo}).
Now we define an axis; 
Choose a start point and an orientation at an edge or a crossing as depicted in Fig. \ref{startpt}. 
Draw a curve according to the following rule ($n=1, 2, 3, \dots $). \\
\noindent $\bullet$ If one enters a $2n$-gon at an edge (resp. crossing), go out at the opposite edge (resp. crossing). (See (a) in Fig. \ref{evenoddregions}.)\\
\noindent $\bullet$ If one enters a $(2n-1)$-gon at an edge (resp. crossing), go out at the opposite crossing (resp. edge). (See (b) in Fig. \ref{evenoddregions}.)\\

\begin{figure}[ht]
\begin{center}
\includegraphics[width=50mm]{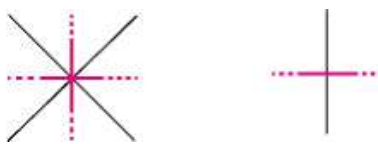}
\caption{Travel transversely. }
\label{pass-ce}
\end{center}
\end{figure}

\noindent $\bullet$ One is permitted to go through each edge transversely at most one time and each crossing transversely at most one time in each direction. (See Fig. \ref{pass-ce}.)\\

Repeat the procedure above until one goes back to the start point with the same direction. 
Thus we obtain a closed curve and call it an {\it axis} of the link projection.

We call each part of an axis bounded by crossings or edges of a link projection a {\it segment}. 
We define the {\it length} of an axis to be the number of segments of the axis. 
For example, in Fig. \ref{ex-axis}, $P$ has four axes with length three, $Q$ has one axis with length 24 and $R$ has axes with length one and three. 
Note that for a knot projection without crossings, the axis is just a point with length zero. 
We have the following proposition: 

\phantom{x}
\begin{proposition}
The total length of all the axes of a link projection $P$ is quadruple of the number of the crossings of $P$. 
\label{fourc}
\end{proposition}
\phantom{x}

\begin{proof}
Let $c$ be the number of crossings of $P$ and $e$ be the number of edges of $P$. 
At each crossing, there are locally four segments of axes of $P$ and at each edge, there are two segments of axes of $P$ as shown in Fig. \ref{pass-ce}. 
The total number of endpoints of all the segments for all the axes is $4c+2e=8c$ (remark that $e=2c$). 
Since each segment has two endpoints, the number of segments is $4c$. 
\end{proof}
\phantom{x}

\noindent We call a segment which lies on a $(2n-1)$-gon (resp. $2n$-gon) an {\it odd-segment} (resp. {\it even-segment}) ($n=1,2,3, \dots $). 
We have the following: 

\phantom{x}
\begin{lemma}
Each axis has an even number of odd-segments. 
\end{lemma}
\phantom{x}

\begin{proof}
Since an axis is a closed curve and each odd-segments has one crossing and one edge on the endpoints, each axis has an even number of odd-segments. 
\end{proof}
\phantom{x}

\noindent Considering all the axes, we have the following: 

\phantom{x}
\begin{corollary}
In total, all the axes go through odd-sided regions even times. 
\label{lem-even}
\end{corollary}
\phantom{x}

\noindent Looking locally, we also have the following (see Fig. \ref{evenoddregions}): 

\phantom{x}
\begin{proposition}
In total, all the axes go through each $n$-gon $n$ times, i.e., 
each $n$-gon of a link projection $P$ has $n$ segments of axes of $P$. 
\label{lem-ngon}
\end{proposition}
\phantom{x}

A link projection $P$ is said to be {\it reducible} if we can put a circle on $S^2$ so that the circle intersects $P$ transversely at just one crossing of $P$. 
We call such a crossing a {\it reducible crossing}. 
A link projection is said to be {\it reduced} otherwise. 
Adams, Shinjo and Tanaka showed in \cite{AST} that every reduced link projection has an even number of odd-sided regions, and Shinjo showed that for both reduced and reducible link projections by giving the equation (1) in \cite{shinjo}. 
We can also see that from Corollary \ref{lem-even} and Proposition \ref{lem-ngon}.

\section{Axis systems}

In this section, we define axis system and prove Theorem \ref{mainthm}. 
Let $P$ be a link projection. 
Give letters to all the regions of $P$, where letters are different each other. 
Let $\alpha$ be an axis of $P$. 
Read up the letters according to $\alpha$ starting from an edge or crossing like Fig. \ref{startpt} with an orientation. 
Give a minus to the head of the sequence we obtained if we started from an edge. 
Thus we obtain a signed sequence, and call it an {\it axis word} of $\alpha$. 
Let $S(P)$ be the set of axis words for all the axes of $P$. 
From Proposition \ref{lem-ngon}, we can see that if a letter $A$ appears $n$ times in $S(P)$, then $A$ corresponds to an $n$-gon of $P$. 
We call a letter an {\it odd-letter} (resp. {\it even-letter}) if the letter corresponds to an odd-sided region (resp. even-sided region). 
Let $\varepsilon \in \{ +, - \}$. 
We call the following permutation $\sigma$ a {\it cyclic permutation on axis words}:
\begin{align*}
\sigma ( \varepsilon A B_1 B_2 \dots B_l ) = 
\begin{cases}
\varepsilon B_1 B_2 \dots B_l A &  \text{(if $A$ is an even-letter)}\\
- \varepsilon B_1 B_2 \dots B_l A &  \text{(if $A$ is an odd-letter)}
\end{cases}
\end{align*}
We call the following permutation $\varphi$ the {\it reversing}: 
$\varphi ( \varepsilon A_1 A_2 \dots A_{l-1} A_l )= \varepsilon A_l A_{l-1} \dots A_2 A_1$. 
Note that $\sigma$ corresponds to changing the start point to read up, and $\varphi$ corresponds to changing the orientation to read up. 
Hence we have the following: 

\phantom{x}
\begin{proposition}
For a link projection $P$, two axis words correspond to the same axis of $P$ if and only if one is obtained from the other one by some $\sigma$ and $\varphi$. 
\end{proposition}
\phantom{x}

We define the {\it axis system} $s(P)$ to be $S(P)$ up to cyclic permutations $\sigma$ and reversing $\varphi$. 
\begin{figure}[ht]
\begin{center}
\includegraphics[width=90mm]{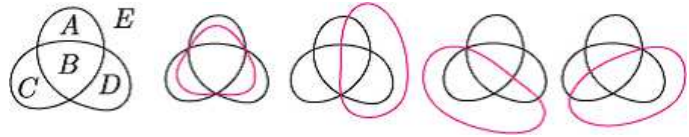}
\caption{Axis system is $\{ ADC, -ABE, BDE, ECB \}$. }
\label{trefoil}
\end{center}
\end{figure}
For example, the knot projection depicted in Fig. \ref{trefoil} has the axis system $\{ ADC, -ABE, BDE, ECB \}$. 
Note that for a knot projection $O$ with no crossings, the axis system is $s(O)= \emptyset$. 
We have the following: 

\phantom{x}
\begin{proposition}
A link projection $P$ is reducible if and only if the axis system $s(P)$ has a word consisting of just one letter or a word including consecutive letters. 
\end{proposition}
\phantom{x}

\begin{figure}[ht]
\begin{center}
\includegraphics[width=30mm]{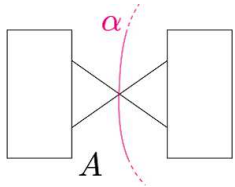}
\caption{An axis $\alpha$ passing a reducible crossing. }
\label{reducible-crossing}
\end{center}
\end{figure}

\begin{proof}
A crossing is a reducible crossing if and only if there are just three regions around it. 
Let $P$ be a reducible link projection. 
Let $\alpha$ be an axis of $P$ which passes a reducible crossing on the same region $A$ as depicted in Fig. \ref{reducible-crossing}. 
If the length of $\alpha$ is one, the word is ``$A$''. 
If the length is two or more, the word includes ``$AA$'' because the axis is divided into two segments at the reducible crossing. 
On the other hand, for reduced link projections, there are no words consisting of just one letter or including consecutive letters because there are four different regions around each crossing, and there are two different regions at the two sides of each edge. 
\end{proof}
\phantom{x}

From now on in this section, we consider only reduced link projections. 
Let $S(P)$ be one of the representations of the axis system $s(P)$ of a link projection $P$. 
As we have seen, from $S(P)$, we can know whether a letter is an odd-letter or an even-letter and whether an axis word begins at a crossing or an edge. 
We give the {\it ce-representation} to each axis word in $S(P)$ as follows: \\

\noindent $\bullet$ If an axis word begins with a minus, give ``$e$'' before the initial letter. 
Otherwise, give ``$c$'' before the initial letter. \\
\noindent $\bullet$ If an odd-letter $A$ has ``$e$'' (resp. ``$c$'') before $A$, give ``$c$'' (resp. ``$e$'') after $A$. \\
\noindent $\bullet$ If an even-letter $A$ has ``$e$'' (resp. ``$c$'') before $A$, give ``$e$'' (resp. ``$c$'') after $A$. \\

\begin{figure}[ht]
\begin{center}
\includegraphics[width=90mm]{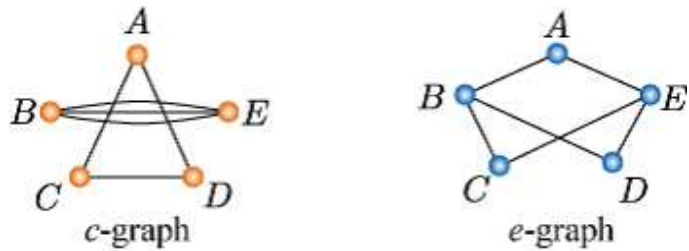}
\caption{The $c$-graph and $e$-graph of $\{ ADC, -ABE, BDE, ECB \}$. }
\label{ce-graph}
\end{center}
\end{figure}

\noindent For example, the ce-representation of $\{ ADC, -ABE, BDE, ECB \}$ is \\
$\{ cAcDcCc, -eAeBcEe, cBeDeEc, cEeCeBc \}$ because $A, C$ and $D$ are even-letters and $B$ and $E$ are odd-letters. 
From ce-representation, we can see the relations between regions. 
For example, if there exists the sequence ``$AcB$'', then the regions $A$ and $B$ share a crossing. 
Similarly, if there exists the sequence ``$AeB$'', then $A$ and $B$ share an edge. 
We represent these relationship visually by {\it c-graph} and {\it e-graph}. 
Let $V= \{ A_1, A_2, \dots , A_l \}$ be the set of all the letters in $S(P)$. 
We assume $V$ as a set of vertices. 
The c-graph (resp. e-graph) is a multi-edge graph such that if two regions corresponding to $A_i$ and $A_j$ share $n$ crossings (resp. $n$ edges), then the corresponding vertices have $n$ edges between them. 
An example is shown in Fig. \ref{ce-graph}. 
Note that $c$-graph and $e$-graph do not depend on the choice of the representation $S(P)$ of $s(P)$. 
Remark that the $c$-graph is coinside with the union of a tait graph without signs of a link projection and the dual of the tait graph (see, for example, \cite{TH}), and the $e$-graph is coinside with the dual of the link projection. 
Hence $c$-graphs and $e$-graphs are planar graphs. 
We have the following:

\phantom{x}
\begin{proposition}
If there exists a crossing of a link projection $P$ which has four regions having the letters $A, B, C$ and $D$ around it with clockwise or counterclockwise rotation, then there exists a cycle with length four $ABCD$ in $e$-graph, and there are edges $AC$ and $BD$ in $c$-graph of $s(P)$.
\label{4cycle}
\end{proposition}
\phantom{x}
\begin{figure}[h]
\begin{center}
\includegraphics[width=90mm]{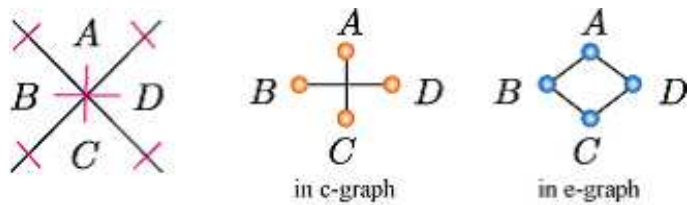}
\caption{If there exists a crossing with regions $A, B, C$ and $D$ around it, then the $e$-graph has the edges $AB, BC, CD$ and $DA$ and the $c$-graph has the edges $AC$ and $BD$. }
\label{ex-crossing}
\end{center}
\end{figure}

\begin{figure}[ht]
\begin{center}
\includegraphics[width=100mm]{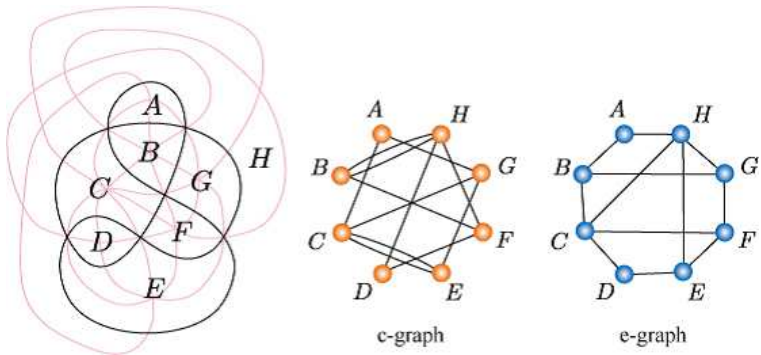}
\caption{A knot projection of the knot $6_3$ and its c-graph and e-graph. }
\label{6-3}
\end{center}
\end{figure}

\noindent (See Fig. \ref{ex-crossing}.) 
Remark that the converse does not hold. 
For example, the knot projection $P$ in the left-hand side in Fig. \ref{6-3} has the axis system 
$$s(P)=\{  AGFDHBCDEGBHFCHGCEFBAHEC \}, $$
and $s(P)$ has the $c$-graph and $e$-graph as depicted in Fig. \ref{6-3}. 
There are nine cycles satisfying the condition of Proposition \ref{4cycle}, and $CBGH, CFGH$ and $CHEF$ do not correspond to a crossing. 
We call such a cycle a {\it dummy}. 
We have the following: 

\phantom{x}
\begin{proposition}
If there are the same number of cycles satisfying the condition of Proposition \ref{4cycle} as the number of letters minus two, we can reobtain the link projection from the $c$-graph and $e$-graph. 
\end{proposition}
\phantom{x}

\begin{proof}
First, the number of letters minus two means the number of crossings of the link projection. 
If there are just the same number of cycles satisfying the condition as the number of crossings of the link projection, there are no dummy cycles. 
Then there is a one-to-one correspondence between the cycles and crossings. 
Since the $e$-graph is coinside with the dual graph of the link projection, each cycle corresponds to a face of the dual graph. 
Then by regarding each cycle as a tile and gluing the tiles on $S^2$, we can obtain the dual graph of the link projection. 
\end{proof}

\begin{figure}[ht]
\begin{center}
\includegraphics[width=100mm]{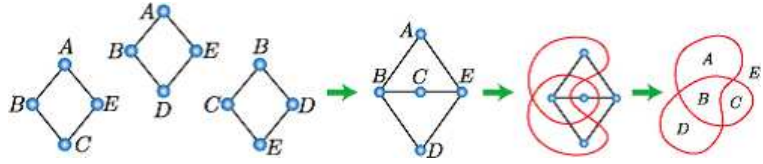}
\caption{We can obtain the knot projection from the $c$-graph and $e$-graph of Fig. 8. }
\label{cycles}
\end{center}
\end{figure}

\noindent See Fig. \ref{cycles} for an example of the method above. 
We prove Theorem \ref{mainthm}. 

\phantom{x}
\noindent {\it Proof of Theorem \ref{mainthm}.} \ 

\begin{figure}[ht]
\begin{center}
\includegraphics[width=100mm]{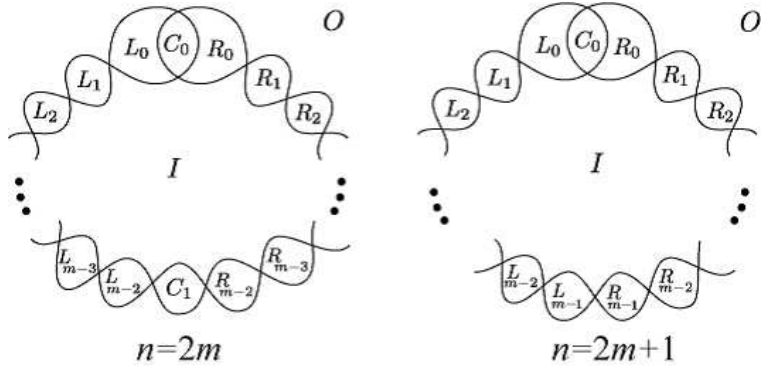}
\caption{The standard projections of twist knots with $n$ crossings. }
\label{twist-eo}
\end{center}
\end{figure}

Give letters to regions of the standard projection of a twist knot with $n$ crossings as shown in Fig. \ref{twist-eo}. 
For the standard projection $P$ of a twist knot with $n=2m$ crossings, we have the following two axes with the words 
$-C_0 R_0 R_1 R_2 \dots R_{m-3} R_{m-2} C_1 L_{m-2} L_{m-3} \dots L_2 L_1 L_0$, $C_0 I C_1 O$, 
and $2(m-2)=n-4$ axes $-R_i I O, \ -L_i I O$ ($i=1,2, \dots , m-2$) and one axis $-L_0 R_0 I O R_0 L_0 O I$. 

For the standard projection $P$ of a twist knot with $n=2m+1$ crossings, we have the two axes 
$-C_0 R_0 R_1 R_2 \dots R_{m-2} R_{m-1} L_{m-1} L_{m-2} \dots L_2 L_1 L_0$, $C_0 I O $, 
and $m-1$ axes $IO, IO, \dots ,IO$ and $m-2$ axes $ -L_i I R_{(m-1)-i} O$ ($i=1,2, \dots , m-2$) and one axis 
$-L_0 R_0 I L_{m-1} O R_0 L_0 O R_{m-1} I$. 
Hence the standard projection $P$ of a twist knot has the axis system $s(P)$ in the statement of Theorem \ref{mainthm}. 
\begin{figure}[ht]
\begin{center}
\includegraphics[width=110mm]{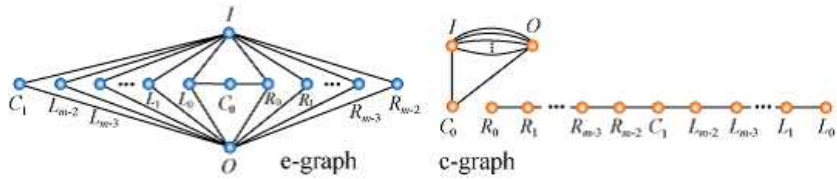}
\caption{The c-graph and e-graph of $s(P)$ for the case of $n=2m$. }
\label{twist-ce}
\end{center}
\end{figure}
\begin{figure}[ht]
\begin{center}
\includegraphics[width=110mm]{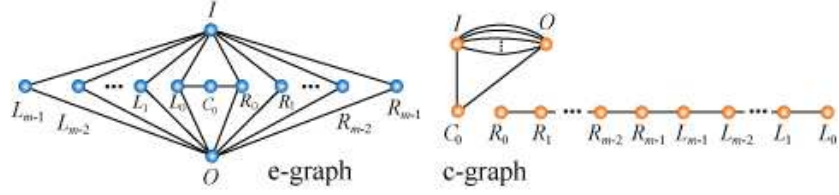}
\caption{The c-graph and e-graph of $s(P)$ for the case of $n=2m+1$. }
\label{twist-ce2}
\end{center}
\end{figure}
On the other hand, there are just $n$ cycles in the e-graph satisfying the condition of Proposition \ref{4cycle}. 
See Fig. \ref{twist-ce} for the case of $n=2m$, and see Fig. \ref{twist-ce2} for the case of $n=2m+1$. 
\hfill$\square$

\phantom{x}

\noindent It is unknown if there exists a one-to-one correspondence between axis systems and link projections. 
If exists, then axis system can be used as a representation of link projections and diagrams such as chord diagram, warping matrix (\cite{shimizu}), warping incidence matrix (\cite{ASW}) and so on. 

\section{Appendix A: Axis polynomials}

\begin{figure}[ht]
\begin{center}
\includegraphics[width=120mm]{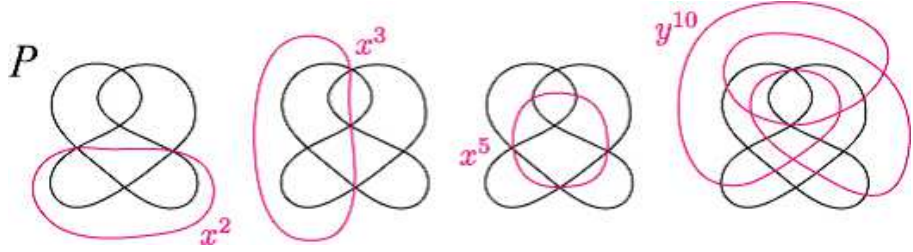}
\caption{The knot projection $S$ has the axis polynomial $A(S; x,y)=x^2+x^3+x^5+y^{10}$. }
\label{twist4}
\end{center}
\end{figure}

In this appendix, we define axis polynomial, which has less information comparing to axis systems, but it is easy to calculate. 
Also, we can derive invariants of link projections from axis polynomials, and as we see in Appendix B, axis polynomial is useful to discuss about symmetry of knot projections. 
We call an axis of a link projection which is a simple closed curve (resp. non-simple closed curve) a {\it simple axis} (resp. {\it non-simple axis}). 
Now we define the {\it axis polynomial}, a two variable polynomial for a link projection $P$ with at least one crossing; 
For a simple axis with length $m$, give $x^m$. 
For a non-simple axis with length $n$, give $y^n$. 
Taking the sum for all the axes, we obtain a polynomial. 
We call it the axis polynomial $A(P; x,y)$ of $P$. 
For example, we have $A(P; x,y)=4x^3$, $A(Q; x,y)=y^{24}$ and $A(R; x,y)=x+x^3$ for the link projections $P$, $Q$ and $R$ in Fig. \ref{ex-axis}. 
For more example, see Fig. \ref{twist4}. 

\noindent By definition, we have the following: 

\phantom{x}
\begin{proposition}
For a link projection $P$, $A(P; 1,1)$ is the number of the axes of $P$. 
\end{proposition}
\phantom{x}

\noindent Similarly, $A(P; 1,0)$ represents the number of the simple axes of $P$. 
Note that axis polynomial, $A(P; 1,1)$ and $A(P; 1,0)$ are invariants of link projections. 
Let $c(P)$ denote the number of crossings of a link projection $P$. 
By Proposition \ref{fourc}, we have the following formula: 

\phantom{x}
\begin{corollary}
$$\frac{\partial}{\partial x}A(P; 1,1)+\frac{\partial}{\partial y}A(P; 1,1)=4c(P).$$
\end{corollary}
\phantom{x}

\noindent From Theorem \ref{mainthm}, we have the following corollary: 

\phantom{x}
\begin{corollary}
Let $P$ be the standard projection of a twist knot with $n$ crossings. 
Then 
$$A(P; x,y)=x^n+x^4+(n-4)x^3+y^8$$
for the case of $n=2m$ and 
$$A(P; x,y)=x^n+x^3+\frac{n-5}{2}x^4+\frac{n-3}{2}x^2+y^{10}$$
for the case of $n=2m+1$ ($m=2,3,4, \dots $). 
\end{corollary}
\phantom{x}

\noindent Note that in \cite{watanabe}, the axis polynomials of the standard projections of torus links are also discussed. 

\phantom{x}
\begin{remark}
There are no one-to-one correspondence between axis polynomials and link projections. 
For example, see Fig. \ref{8crossings}. 
\begin{figure}[ht]
\begin{center}
\includegraphics[width=110mm]{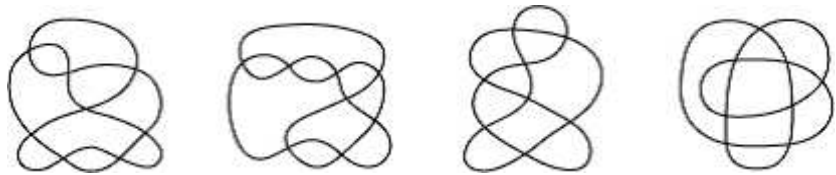}
\caption{Projections of the knots $8_7$, $8_9$, $8_{13}$ and $8_{17}$ have the same axis polynomial, $y^{32}$. }
\label{8crossings}
\end{center}
\end{figure}

\end{remark}

\section{Appendix B: Symmetry}

Axes are useful to consider the symmetry of knot projections. 
A link projection $P$ is {\it symmetric with respect to a circle $\mu$} on $S^2$ (or just {\it symmetric}) if there exists a homeomorphism $\varphi : S_+^2 \to S_-^2$ such that $\varphi (x)=x$ for any $x \in \mu$ and $\varphi (P \cap S_+^2 )=P \cap S_-^2$, where $S_+^2$ and $S_-^2$ are 2-discs such that $S_+^2 \cup S_-^2 = S^2$ and $S_+^2 \cap S_-^2 = \mu$. 
We have the following: 

\phantom{x}
\begin{proposition}
If a link projection $P$ is symmetric with respect to a circle $\mu$ on $S^2$, the circle $\mu$ is corresponding to a component of the link or a simple axis. 
In particular, if a knot projection $P$ with at least one crossing is symmetric with respect to a circle $\mu$, then $\mu$ is a simple axis of $P$. 
\label{prop31}
\end{proposition}
\phantom{x}

\begin{figure}[ht]
\begin{center}
\includegraphics[width=40mm]{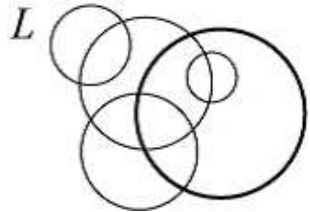}
\caption{The link projection $L$ is symmetric with respect to the thick component itself. }
\label{symmetric-link}
\end{center}
\end{figure}

\noindent (See Fig. \ref{symmetric-link} for a link projection.) 
Remark that the converse of Proposition \ref{prop31} does not hold. 
For example, the knot projection $S$ in Fig. \ref{twist4} has three simple axes; $A(S; 1,0)=3$. 
However, $P$ is not symmetric with respect to the simple axis with length two. 
We can see that by counting the number of 2-gons at each side of the axis. 
More simply, we can say the following: 

\phantom{x}
\begin{corollary}
If a knot projection $P$ is symmetric, $P$ has a simple axis. 
\label{symmetric-simple}
\end{corollary}
\phantom{x}

\noindent The contraposition of Corollary \ref{symmetric-simple} is useful. 
For example, we can say that the knot projection $Q$ in Fig. \ref{ex-axis} is not symmetric because $A(Q; 1,0)=0$.

\section*{Acknowledgments}
The authors are grateful to Taira Akiyama, Yukari Fuseda, Naoya Irisawa, Yui Onoda and Sayaka Shimizu for valuable discussions and conversations at seminars in NIT, Gunma College. 
The second author thanks Professors Makiko Ishiwata, Takako Kodate and Kazuaki Kobayashi for helpful discussions and advice on knot projections. 
She was partially supported by Grant for Basic Science Research Projects from the Sumitomo Foundation. 
Finally, the authors are deeply grateful to the referee for helpful information and suggestions.

\end{document}